\def\timestamp{%
Time-stamp: <betaN-problems.tex: donderdag 17-10-2024 at 15:52:21 (cest)>}
\def\stripname Time-stamp: <#1: #2 #3>{#3}
\edef\filedate{\expandafter\stripname\timestamp}

\documentclass[a4paper]{amsart}
\usepackage[margin=30mm]{geometry}
\newcounter{question}

\newenvironment{question}%
   {\par\medbreak\refstepcounter{question}
    \noindent\llap{\textbf{\thequestion}. }\ignorespaces}%
   {\endgraf
    \medbreak}
\newcommand\comments
   {\endgraf
    \everypar{\setbox0=\lastbox\everypar{}\textsc{Comments}: }}

\DeclareMathSymbol\C0{AMSb}{`C}
\DeclareMathSymbol\Gap0{AMSb}{`G}
\DeclareMathSymbol\HH0{AMSb}{`H}
\DeclareMathSymbol\I0{AMSb}{`I}
\DeclareMathSymbol\N0{AMSb}{`N}
\DeclareMathSymbol\Q0{AMSb}{`Q}
\DeclareMathSymbol\R0{AMSb}{`R}
\DeclareMathSymbol\Z0{AMSb}{`Z}
\DeclareMathSymbol\diamond0{AMSa}{"06}
\newcommand\betaN{\beta\N}
\newcommand\Nstar{\N^*}
\newcommand\Hstar{\HH^*}

\DeclareMathSymbol\restr\mathbin{AMSa}{"16}

\newcommand\RPC{\operatorname{RPC}}

\newcommand\Triv{\mathsf{Triv}}
\newcommand\Aut{\mathsf{Aut}}
\newcommand\RK{\mathbf{RK}}
\newcommand\RF{\mathbf{RF}}
\newcommand\Cantor[1][\omega]{\vphantom{2}^{#1}2}
\newcommand\Cantorc{\Cantor[\cee]}
\newcommand\cf{\operatorname{cf}}
\newcommand\Fn{\operatorname{Fn}}
\newcommand\id{\operatorname{Id}}
\newcommand\pow{\mathcal{P}}
\newcommand\fin{\mathrm{fin}}
\newcommand\powNfin{\pow(\N)/\fin}
\newcommand\calA{\mathcal{A}}
\newcommand\calB{\mathcal{B}}
\newcommand\calF{\mathcal{F}}
\newcommand\calI{\mathcal{I}}
\newcommand\calO{\mathcal{O}}
\newcommand\calU{\mathcal{U}}
\newcommand\calX{\mathcal{X}}
\newcommand\calY{\mathcal{Y}}

\newcommand\axiom{\mathsf}
\newcommand\CH{\axiom{CH}}
\newcommand\MA{\axiom{MA}}
\newcommand\NCF{\axiom{NCF}}
\newcommand\PFA{\axiom{PFA}}
\newcommand\ZFC{\axiom{ZFC}}
\newcommand\ZF{\axiom{ZF}}
\newcommand\OCA{\axiom{OCA}}
\newcommand\MAnotCH{{\MA+\neg\CH}}

\newcommand\bee{\mathfrak{b}}
\newcommand\cee{\mathfrak{c}}
\newcommand\dee{\mathfrak{d}}
\newcommand\enn{\mathfrak{n}}
\newcommand\you{\mathfrak{u}}

\newcommand\cl{\operatorname{cl}}
\newcommand\clos{\overline}
\newcommand\card[1]{\mathopen|#1\mathclose|}
\let\abs=\card

\newcommand\orpr[2]{\langle{#1},{#2}\rangle}
\newcommand\RKeq{\mathrel{\equiv_\RK}}
\newcommand\leRK{\mathrel{\le_\RK}}
\newcommand\leRF{\mathrel{\le_\RF}}

\newcommand\lRF{\mathrel{<_\RF}}

\newcommand\preim{^\gets}
\newcommand\Spchi{\mathrm{Sp}_\chi}

\newcommand\pri[1]{}

\DeclareMathSymbol\streepje0{operators}{"2D}
\newcommand\plim[1][p]{#1\streepje\!\lim}

\usepackage{amsrefs}

\begin{document}

\title{Problems on $\betaN$}

\author{Klaas Pieter Hart}

\address{Faculty EEMCS\\TU Delft\\
         Postbus 5031\\2600~GA {} Delft\\the Netherlands}
\email{k.p.hart@tudelft.nl}
\urladdr{http://fa.ewi.tudelft.nl/\~{}hart}

\author{Jan van Mill}
\address{KdV Institute for Mathematics\\
         University of Amsterdam\\
         P.O. Box 94248\\
         1090~GE {} Amsterdam\\
         The Netherlands}
\email{j.vanmill@uva.nl}
\urladdr{https://staff.fnwi.uva.nl/j.vanmill}

\date{\filedate}

\begin{abstract}
This is an update on, and expansion of, our paper
\textsl{Open problems on $\beta\omega$}
in the book \textsl{Open Problems in Topology}.
\end{abstract}


\maketitle

\section*{Introduction}

In 1990 we contributed a paper to the book
\textsl{Open problems in Topology}, \cite{MR1078636}, titled
\textsl{Open problems on~$\beta\omega$} (\cite{MR1078643}).

Through the years some of these problems were solved, some were shown to be
related to other problems, and some are still unsolved.
In the first years after the publication of the book there were regular
updates on the problems in the journal
\textsl{Topology and its Applications};
in~2004 these were collated and extended in a comprehensive status report,
\cite{MR2023411},
by Elliott Pearl.

The COVID-19 pandemic provided a good opportunity to go through our original
paper again and provide a new update of the status of the problems as well as
to collect and formulate new questions on the fascinating object that
is~$\betaN$.

Many of the comments below incorporate information from Elliott Pearl's
update mentioned above, but there have, of course, been many developments
in the years since.

The numbering of the problems is different from that in the first paper
because we have moved some questions around and combined related questions
into more comprehensive problems.
We have made no attempt to separate the solved problems from the unsolved ones.
We wanted to keep related problems together and even though we consider a
problem solved the reader may disagree and be inspired to investigate
variations or strengthenings of the answers.

We should also mention the book
\textsl{Open problems in Topology II}~\cite{MR2367385},
edited by Elliot Pearl, that contains many more problems in topology,
and in particular a paper by Peter Nyikos,
\textsl{\v{C}ech-Stone remainders of discrete spaces},
that, as the title indicates, deals with problems on~$\beta\kappa$ for
arbitrary infinite cardinals~$\kappa$.

\section{Preliminaries}

The main objects of study in this paper are the space~$\betaN$ and its
subspace~$\Nstar$.

For a quick overview of their properties we refer to Chapter~D-18
of~\cite{MR2049453}; the paper~\cite{MR776630} offers a more comprehensive
treatment.
Here we collect some of the basic facts about our spaces in order to fix
some notation that will be used throughout the paper.

To begin: $\betaN$ is the set of ultrafilters on the set~$\N$ of natural
numbers, endowed with the topology generated by the base
$\{\clos{A}:A\subseteq\N\}$, where $\clos{A}$~denotes the set of ultrafilters
that contain~$A$.
The readily established equality
$\betaN\setminus\clos{A}=\clos{\N\setminus A}$
confirms what the notation~$\clos{A}$ suggests:
the set $\clos{A}$ is open \emph{and} closed, and also equal to the closure
of~$A$ in~$\betaN$.

We identify an element~$n$ of~$\N$ with the ultrafilter $\{A:n\in A\}$ and
thus consider $\N$ to be a subset of~$\betaN$.
The complement $\betaN\setminus\N$ is the set of \emph{free} ultrafilters
on~$\N$ and is denoted~$\Nstar$; we extend this notation to all subsets
of~$\N$ and write $A^*=\clos{A}\setminus A$ whenever $A\subseteq\N$.

A map~$\varphi$ from~$\N$ to itself induces a map~$\beta\varphi$
from~$\betaN$ to itself: $\beta\varphi(u)$~is the ultrafilter
generated by~$\{\varphi[A]:A\in u\}$.

\section{Autohomeomorphisms}

The autohomeomorphisms of~$\betaN$ correspond to the permutations of~$\N$ and
are, as such, not very interesting topologically.
The autohomeomorphisms of~$\Nstar$ offer more challenges.

In what follows $\Aut$ denotes the autohomeomorphism group of~$\Nstar$,
and $\Triv$~denotes the subgroup of trivial autohomeomorphisms.
Here a \emph{trivial} autohomeomorphism is one with an `easy' description:
an autohomeomorphism~$h$ of~$\Nstar$ is trivial if there are co-finite
subsets~$A$ and~$B$ of~$\N$ and a bijection~$\varphi:A\to B$ such that
$h=\varphi^*$, where $\varphi^*$~denotes the restriction of~$\beta\varphi$
to~$\Nstar$.

\begin{question}\label{q.triv.normal}
Can $\Triv$  be a proper normal subgroup of $\Aut$, and if \emph{yes} what is
(or can be) the structure of the factor~group $\Aut/\Triv$;
and if \emph{no} what is (or can be) $[\Triv:\Aut]$?
\comments
A very concrete first step would be to investigate what can one say about an
autohomeomorphism~$h$ that satisfies $h^{-1}\circ\Triv\circ h=\Triv$.

A related question: what is the minimum number of autohomeomorphisms
necessary to add to $\Triv$ to get a generating set for~$\Aut$?

Of course that number is~$0$ when $\Aut=\Triv$,
but can it be non-zero and finite?
\end{question}

If $h\in\Aut$ then $I(h)$ denotes the family of subsets of~$\omega$
on which $h$~is trivial, that is, $A\in I(h)$ iff there is a function
$h':A\to\omega$ such that $h(B^*)=h'[B]^*$ whenever $B\subseteq A$.

If $I(h)$ contains an infinite set then $h$~is \emph{somewhere trivial},
otherwise $h$~is totally non-trivial.

The ideal $I(h)$ determines an open set $O_h$: the union
$\bigcup\{A^*:A\in I(h)\}$; its complement~$F_h$ is closed and could be
said to be the set of points of~$\Nstar$ where $h$~is truly
non-trivial.

\begin{question}
Does the existence of a (totally) non-trivial automorphism imply that
$\Aut$ is simple?
\comments
This question asks more than the opposite of question~\ref{q.triv.normal};
a yes answer here would imply a no answer there, but a negative answer
there could go together with a negative answer here.
\end{question}

\begin{question}
Is it consistent with $\MA+\neg\CH$ that a totally non-trivial
automorphism exists?
\comments
The answer yes.
This was established by Shelah and Stepr\=ans in~\cite{MR1896046}.
Consistency is the best one can hope for: 
in~\cite{MR935111} Shelah and Stepr\=ans proved that $\PFA$~implies
all autohomeomorphisms of~$\Nstar$ are trivial; they also indicated how
the implicit large cardinal assumption can be avoided and 
use $\diamond$ on~$\omega_2$ to capture and eliminate any potential non-trivial
autohomeomorphisms in a countable support iteration of length~$\omega_2$. 
Though not stated explicitly by the authors it is clear that one can modify
the iteration so as to obtain a model that satisfies $\MA_{\aleph_1}$ as well.
In~\cite{MR1202874} Veli\v{c}kovi\'c showed that the conjunction 
of~$\MA_{\aleph_1}$ and~$\OCA$ implies that all autohomeomorphisms are trivial.
\end{question}

\begin{question}
Is it consistent to have a non-trivial automorphism, while for every
$h\in\Aut$ the ideal~$I(h)$ is the intersection of finitely many prime ideals?
\comments
In topological terms: can one have non-trivial autohomeomorphisms but only
very mild ones; the set of points where an autohomeomorphism is truly
non-trivial is always finite.
\end{question}

\begin{question}
Is every ideal $I(h)$ a $P$-ideal? 
\comments
This was asked explicitly in~\cite{MR935111}*{Question~2} in case every
autohomeomorphism is somewhere trivial, after it was shown that $\PFA$
implies a yes answer.
However, as mentioned above, $\PFA$ implies that all autohomeomorphisms 
are trivial, so that $I(h)$~is, in fact, always an improper ideal.

Of course this question only makes sense in case $I(h)$ is not equal
to the ideal of finite sets.
Also, if \emph{every} autohomeomorphism is somewhere trivial then
every $I(h)$~is a tall ideal and hence the set of points of non-triviality
is nowhere dense.

Without the additional condition that every autohomeomorphism is somewhere
trivial the answer is consistently negative.
The Continuum Hypothesis lets one construct an autohomeomorphism~$h$ that is 
trivial, in fact the identity, on the members of a partition~$\calA$ 
of~$\omega$ into infinite sets, and so that there is a point~$u$ on the 
boundary of $\bigcup\{A^*:A\in\calA\}$ such that $h$~is not trivial on 
each neighbourhood of~$u$. 
This implies there is no~$B\in I(h)$ such that $A\subseteq B$ 
for all~$A\in\calA$.
\end{question}

\begin{question}
If every automorphism is somewhere trivial, is then every automorphism
trivial?
\comments
This is undecidable.

Shelah proved the consistency of ``all autohomeomorphisms are trivial'',
see~\cite{MR1623206}.
Shelah and Stepr\=ans proved the consistency with $\MA_{\aleph_1}$
of ``every autohomeomorphism is somewhere trivial, yet there is a non-trivial
autohomeomorphism'' in~\cite{MR1271551};
as noted above they proved in~\cite{MR1896046} that $\MA$~does not imply
that all autohomeomorphisms are somewhere trivial.
\end{question}

Given a cardinal $\kappa$ call an autohomeomorphism~$h$
\emph{weakly $\kappa$-trivial} if the set
$\{p:p\not\RKeq h(p)\}$ has cardinality less than~$\kappa$.
Here $p\RKeq q$~means that $p$ and~$q$ have the same type, i.e.,
$q=\pi^*(p)$ for some permutation~$\pi$ of~$\N$.

\begin{question}
For what cardinals $\kappa$ is it consistent to have that all
autohomeomorphisms are weakly $\kappa$-trivial?
\comments
Since a trivial autohomeomorphism is weakly $1$-trivial we see
that $\kappa=1$ is a possibility.
And of course the candidates are less than or equal to~$2^\cee$.
\end{question}

\begin{question}
If $h$ is weakly $1$-trivial is $h$ then trivial?
\comments
This is a uniformization question: if for every $p\in\Nstar$
there is a permutation~$\pi_p$ such that $h(p)=\pi_p^*(p)$ is there then
one (almost) permutation~$\pi$ of~$\Nstar$ such that $h(p)=\pi^*(p)$ for
all~$p$?
\end{question}

\begin{question}
($\MA+\lnot\CH$)
if $p$ and $q$ are $P_\cee$-points\pri{Pc-point@$P_\cee$-point}
is there an $h$ in $\Aut$ such that $h(p)=q$?
\comments
This is undecidable.

Shelah and Stepr\=ans proved the consistency of $\MA+\lnot\CH$ with
``all autohomeomorphisms are trivial''
in~\cite{MR935111}; in this model there are $\cee$~many autohomeomorphisms
and $2^\cee$~many $P_\cee$-points.

Stepr\=ans proved the consistency of a positive answer in~\cite{MR1239060}.
\end{question}

In the investigations into the previous question the following
equivalence relation was used: $p\equiv q$ means that there are two partitions
$\{A_n:n\in\omega\}$ and~$\{B_n:n\in\omega\}$ of~$\N$ into finite sets such that
$$
(\forall P\in p)(\exists Q\in q)(\forall n)
\bigl(\card{P\cap A_n}=\card{Q\cap B_n}\bigr)
$$
The following question was left open.

\medskip

\begin{question}
Is $\equiv$ different from $\RKeq$ in $\ZFC$\pri{ZFC@\ZFC}?
\comments
\end{question}

\begin{question}
Are the autohomeomorphisms of $\Nstar$ induced by the shift map
$\sigma:n\mapsto n+1$ and by its inverse conjugate?
\comments
Recently Will Brian showed that the answer to this question is affirmative
assuming~$\CH$, see~\cite{arxiv:2402.04358}.

See~\cites{MR1035463,MR4138425} for earlier results.
Under~$\CH$ the autohomeomorphism group of $\Nstar$ is simple, yet it has
the maximum possible number of conjugacy classes:~$2^\cee$.
This suggests questions about the number and nature of conjugacy classes
of this group, in $\ZFC$ or under various familiar extra set-theoretical
assumptions, see also~\cite{MR4560745}.
\end{question}

\begin{question}
Does $\Nstar$ have a universal autohomeomorphism?
\comments
This is a question with many possible variations.
The definition of universality that we adopt here is as follows:
$f:\Nstar\to\Nstar$ is universal if for every closed subspace~$F$ of~$\Nstar$
and every autohomeomorphism~$g$ of~$F$ there is an embedding $e:F\to\Nstar$
such that $g=e^{-1}\circ (f\restr e[F])\circ e$.

One can ask whether there is a universal autohomeomorphism at all,
whether the shift~$\sigma$ is universal (for autohomeomorphisms without
fixed points),
whether there is a universal autohomeomorphism just for autohomeomorphisms
without fixed points.

The authors have shown that there is a universal autohomeomorphism of~$\Nstar$
under~$\CH$ and that there is no trivial universal autohomeomorphism.
See~\cite{MR4398473}.

We also note that $\N$ has a universal permutation:
take a permutation of~$\N$ that has infinitely many $n$-cycles, for every~$n$,
and infinitely many infinite cycles (copies of~$\Z$ with the shift).
Every other permutation of~$\N$ can be embedded into this one.
\end{question}

\section{Subspaces}

\begin{question}
For what $p$ are $\Nstar\setminus\{p\}$
and $\betaN\setminus\{p\}$ non-normal?
\comments
Originally this question had the word `equivalently' after the `and' 
(in parentheses).
Since $\Nstar\setminus\{p\}$ is closed in $\betaN\setminus\{p\}$ there is an
implication between the non-normality of these spaces but we do not know
whether that implication is reversible.
Thus this question may actually be two separate ones.

Under $\CH$ the answer is, in both cases, ``for every point'',
see \cites{MR234422,MR321012,MR292035,MR861501}.
There are some results for some special types of points,
see, e.g., B\l aszczyk and Szyma\'nski~\cite{BlaszczykSzymanski1980},
Gryzlov~\cite{MR760274}, and Logunev~\cite{MR4405820}, but
a general answer is wanting.
\end{question}

\begin{question} Is it consistent that there is a non-butterfly
          point\pri{non-butterfly point} in $\Nstar$?
\comments
We call $p$ a butterfly point if there are disjoint sets~$A$ and~$B$ such
that $p$~is the only common accumulation point of~$A$ and~$B$,
that is: $A^d\cap B^d=\{p\}$.

The points used by B\l aszczyk and Szyma\'nski~\cite{BlaszczykSzymanski1980}
in their (partial) answer to the previous question are easy-to-describe
butterfly points.
Let $X=\{x_n:n\in\omega\}$ be a discrete subset of~$\Nstar$,
let $A=\cl X\setminus X$ and take $p\in A$.
Let $\{B_n:n\in\omega\}$ be a partition of~$\N$ such that $B_n\in x_n$ for
all~$n$ and let $q$ be the ultrafilter $\{Q:p\in\cl\{x_n:n\in Q\}\}$.
The set $B=\bigcap_{Q\in q}\bigl(\bigcup\{B_n:n\in Q\}\bigr)^*$ is closed
and $\{q\}=A^d\cap B^d$.
Thus butterfly points exist.

By contrast, in \cite{MR1058805} Be\v{s}lagi\'c and Van Douwen showed that
it is consistent with all consistent cardinal arithmetic that all points
of~$\Nstar$ are butterfly points.
\end{question}

\begin{question}
Is it consistent that $\Nstar\setminus\{p\}$ is $C^*$-embedded in~$\Nstar$
for some but not all~$p\in\Nstar$?
\comments
The answer is yes: in~\cite{MR1434375} Alan Dow showed that
in the Miller model $\Nstar\setminus\{p\}$ is $C^*$-embedded iff $p$~is not
a $P$-point.
There are $P$-points in the Miller model: every ground-model $P$-point
generates a $P$-point in the extension.
\end{question}

\begin{question}
What spaces\pri{subspace!of $\beta\omega$} can be embedded in $\beta\omega$?
\comments
This is a very general question and a definitive answer looks out of reach
for now, even for closed subspaces.

The Continuum Hypothesis implies that the closed subspaces of $\beta\omega$
are exactly the compact zero-dimensional $F$-spaces; in fact, these are also
exactly the closed $P$-sets in~$\Nstar$.
The implication does not reverse: in~\cite{MR1152978} it is shown that
every compact zero-dimensional $F$-space is a (closed) subspace of~$\Nstar$
in any model obtained by adding $\aleph_2$~many Cohen reals to a
model of~$\CH$.

Dow and Vermeer proved in~\cite{MR1137221} that it is consistent that
the $\sigma$-algebra of Borel sets of the unit interval is not the quotient
of any complete Boolean algebra.
By Stone duality, this yields a compact basically disconnected space,
hence a compact zero-dimensional $F$-space, of weight~$\cee$ that cannot
be embedded into any extremally disconnected space, in particular
not into~$\betaN$.

Some $\ZFC$ results are available.
For instance: if $X$~is a compact space of countable cellularity
that is a continuous image of~$\Nstar$ then its projective cover~$E(X)$
can be embedded in~$\Nstar$ as a $\cee$-OK set (a weakening of the notion
of a $P$-set).
This was proved by van~Mill in~\cite{MR637426} and applies to all separable
compact extremally disconnected spaces as well as to the projective covers
of Suslin lines and of Bell's ccc non-separable remainder~\cite{MR624458}.

Van Douwen proved in unpublished work that every $P$-space of weight~$\cee$
(or less) can be embedded into~$\betaN$.
In fact he proved that for every infinite cardinal~$\kappa$ every $P$-space
of weight~$2^\kappa$ can be embedded in~$\beta\kappa$.
The argument was sketched and extended in~\cite{MR674103} and we summarize
it here for the reader's convenience.

Let $X$ be a $P$-space of weight~$2^\kappa$ and embed it into the Cantor cube
$C=2^{2^\kappa}$ of weight~$2^\kappa$.
Next consider the projective cover~$\pi:E(C)\to C$ of this cube.
The Cantor cube is a group under coordinatewise addition modulo~$2$, so
for every $p\in C$ the map $\lambda_p:x\mapsto x+p$ is a homeomorphism; this
homeomorphism lifts to a homeomorphism $\Lambda_p:E(C)\to E(C)$ with the
property that $\pi\circ\Lambda_p=\lambda_p\circ\pi$.
Now take one point~$u_0\in E(C)$ that maps to the neutral element~$0$ of~$C$
and consider the subspace $X'=\{\Lambda_p(u_0):p\in X\}$ of~$E(C)$.
Using the fact that regular open sets in~$C$ are, up to permutation of the
coordinates, of the form $U\times2^I$ where $U$~is regular open
in the Cantor set~$2^\omega$ and $I=2^\kappa\setminus\omega$, one shows
that $\pi$~is actually a homeomorphism from~$X'$ to~$X$.
Finally then, as $\pi$~is irreducible and $C$~has density~$\kappa$,
the density of~$E(C)$ is equal to~$\kappa$ as well.
Therefore there is a continuous surjection $f:\beta\kappa\to E(C)$ and one
can take a closed subset~$F$ of~$\beta\kappa$ such that $f\restr F$~is
irreducible and onto.
As $E(C)$~is extremally disconnected this restriction is a homeomorphism
and we find our copy of~$X$ in~$F$.

The extension in~\cite{MR674103} delivers more but at a cost:
one embeds $\beta X$ in a suitable Cantor cube, possibly of a larger
weight than that of~$X$ itself.
What this delivers is that the copy of~$X$ in $\beta\lambda$
(where $\lambda$ may be larger than the $\kappa$ above) is $C^*$-embedded.

Thus we get the general statement that every $P$-space can be $C^*$-embedded
in a compact extremally disconnected space.

This argument also shows that $2^{\aleph_0}=2^{\aleph_1}$ implies that
$\beta\omega_1$ embeds into~$\betaN$.
For $\beta\omega_1$ embeds into the Cantor cube $2^{2^{\omega_1}}$, which under
our assumption is the same as~$2^\cee$.
The latter is a continuous image of~$\betaN$ and an irreducible preimage
of~$\beta\omega_1$ will be homeomorphic to~$\beta\omega_1$.
If $2^{\aleph_0}<2^{\aleph_1}$ then $\beta\omega_1$ can not be embedded
into~$\Nstar$ because its weight, which is $2^{\aleph_1}$, is larger than that
of~$\Nstar$.
\end{question}



\begin{question}
Describe the closed $P$-sets\pri{P-set@$P$-set} of $\Nstar$. \label{q.P-sets}
\comments
This has a quite definitive answer under $\CH$: every compact zero-dimensional
$F$-space of weight $\mathfrak{c}$ can be embedded in~$\Nstar$ as a $P$-set.
What we are looking for are properties that can be established in~$\ZFC$,
or provably can not.
For example: one cannot prove in~$\ZFC$ that there is a $P$-set homeomorphic
to~$\Nstar$ itself, see~\cite{MR976360}, or that there is a $P$-set
that satisfies the ccc, see~\cite{MR1253914}.

One can ask if cellularity less than~$\cee$ is at all possible.

There are various nowhere dense closed $P$-sets that one can write down
explicitly.
To give two familiar examples, among many, we consider the density
ideal~$\calI_d$ and the summable ideal~$\calI_\Sigma$.
The first is defined by
$$
I\in\calI_d \text{ \quad iff \quad}\lim_{n\to\infty}\frac1n|A\cap n|=0
$$
and
the second as
$$
I\in\calI_\Sigma \text{ \quad iff \quad}\sum_{n\in A}\frac1n \text{ converges.}
$$
These ideals have been studied widely but we would like to know: what are
the topological properties of the nowhere dense closed $P$-sets
$$
F_d=\Nstar\setminus\bigcup\{A^*:A\in\calI_d\}
\text{ \quad and \quad}
F_\Sigma=\Nstar\setminus\bigcup\{A^*:A\in\calI_\Sigma\}
$$
Rudin established in~\cite{MR0216451}
that $F_d$ contains no $P$-points and even that no countable
set of $P$-points accumulates at a point of~$F_d$.
Indeed let $u$ be a $P$-point and observe first that for every~$n$ there
is an $i_n<n$ such that $U_n=\{m\in\N: m\equiv i_n\pmod{n}\}$ belongs
to~$u$.
Because $u$~is a $P$-point there is then a $U\in u$ such that $U\subseteq^*U_n$
for all~$n$.
But this implies $U\in\calI_d$, so $u\notin F_d$.
Because $F_d$ is a $P$-set this implies that no countable set of $P$-points
has accumulation points in~$F_d$.

There are certain similarities between the two sets and $\N^*$ itself.
Consider the map $f:\N\to\N$, defined by $f(n)=k$ iff $k!<n\le(k+1)!$.
It is an elementary exercise to show that
$$
\limsup_{n\to\infty}\frac1n|f\preim[X]\cap n|=1
\text{ \quad en \quad}
\sum_{n\in f\preim[X]}\frac1n=\infty
$$
whenever $X$~is an infinite subset of~$\N$.
This implies that $\beta f$ maps both $F_d$ and $F_\Sigma$ onto~$\Nstar$
and it allows for the lifting of many combinatorial structures on~$\Nstar$
to these sets.
It is clear that the restriction of~$\beta f$ to $\Nstar$ is an open map
onto~$\Nstar$ itself, whether its restrictions to $F_d$ and~$F_\Sigma$
are open as well is less clear.
\end{question}

\begin{question}
Which compact zero-dimensional $F$-spaces admit an open map onto~$\Nstar$?
\comments
This question is related to Van Douwen's paper~\cite{MR1062775}, where open
maps are used to transfer information from $\Nstar$ to other remainders.
As a special case one can investigate whether the sets~$F_d$ and~$F_\Sigma$
from the Question~\ref{q.P-sets} admit open maps onto~$\Nstar$
(if the map~$\beta f$ given there does not already give open maps).
\end{question}

\begin{question}
Is there a nowhere dense copy of $\Nstar$ in $\Nstar$ that is
a $\cee$-OK-set in $\Nstar$?\label{q.c-OK}
\comments
Alan Dow showed in~\cite{MR3209343} that there a nowhere dense copy
of~$\Nstar$ that is not of the form~$\cl{D}\setminus D$ for some countable
and discrete subset~$D$ of~$\Nstar$.
This was later improved by Dow and van Mill in~\cite{MR4518082} to a nowhere
dense copy that is a weak $P$-set.
In light of the comments for question~\ref{q.P-sets} the present question
asks for the best that we can get in~$\ZFC$.
Most likely the answer to this question will require a new idea as the
constructions in the papers cited above produce sets that are definitely
not $\cee$-OK in~$\Nstar$.
\end{question}

\begin{question} Is every subspace of $\Nstar$ strongly zero-dimensional%
\pri{strongly zero-dimensional|see{zero-dimensional, strongly}}%
\pri{zero dimensional!strongly}?
\comments
It is clear that every subspace is zero-dimensional and that closed subspaces
are even strongly zero-dimensional, but for general subspaces this question
is quite open.
Until recently it was not even known whether there was an example of a
zero-dimensional $F$-space that is not strongly zero-dimensional,
see~\cite{MR4384168}.

If the answer is negative then a secondary question suggests itself immediately:
is there an upper bound to the covering dimension of subspaces of~$\Nstar$?
\end{question}

\begin{question}
Is every nowhere dense subset of $\Nstar$%
\pri{subset!nowhere dense!of $\Nstar$}%
\pri{nowhere dense subset|see{subset, nowhere dense}} a
$\cee$-set\pri{c-set@$\cee$-set}?\label{q.cee-set}
\comments
In general a set~$A$ is called a $\kappa$-set if there is a pairwise disjoint
family~$\calO$ of open sets of cardinality~$\kappa$ and
such that $A\subseteq\bigcap\{\cl O:O\in\calO\}$.

That the answer is positive is called by some ``The $\cee$-set conjecture''.
In~\cite{MR1056181} Simon proved that this question is the same as
``Is there a maximal nowhere dense subset in~$\Nstar$?''.
The questions are the same in that the answer ``no'' to one is equivalent to
the answer ``yes'' to the other:
Every nowhere dense set in $\Nstar$ is a $\cee$-set if
and only if every nowhere dense set in $\Nstar$ is a
nowhere dense subset of another nowhere dense set
(this is the order that we are considering).

There is a purely combinatorial reformulation of this question,
denoted $\RPC(\omega)$ in~\cite{MR991597}: if $\calA$ is an infinite maximal
almost disjoint family then $\calI^+(\calA)$ has an almost disjoint refinement.
Here, $\calI^+(\calA)$ is the family of sets not in the ideal~$\calI(\calA)$
generated by~$\calA$ and the finite sets and an almost disjoint refinement
is an almost disjoint family~$\calB$ with a map $X\mapsto B_X$
from $\calI^+(\calA)$ to~$\calB$ such that $B_X\subseteq^* X$ for all~$X$.

Finally, we should mention that the answer is positive for one-point sets:
all points of $\Nstar$ are $\cee$-points, see~\cite{MR567994}.
\end{question}

\begin{question}
Does there exist a completely separable maximal almost disjoint family?
\comments
This question is related to Question~\ref{q.cee-set} because by
\cite{MR991597}*{Theorem~4.19} a positive answer to that question
implies the existence of an abundance of completely separable maximal almost
disjoint families; where a maximal almost disjoint family~$\calA$
is \emph{completely separable} if it is itself an almost disjoint
refinement of~$\calI^+(\calA)$.

Whether completely separable maximal almost disjoint families exist is a
problem first raised by Erd\H{o}s and Shelah in~\cite{MR319770}.

Currently the best result is due to Shelah who showed
in~\cite{MR2894445} that the answer is positive if $\cee<\aleph_\omega$
and that a negative solution would imply consistency of the existence of
large cardinals.

It is not (yet) clear whether this question and Question~\ref{q.cee-set}
are equivalent.
Thus far \emph{constructions} of completely separable maximal almost disjoint
families (in some model or another) could always be adapted to
prove~$\RPC(\omega)$, but there is currently no proof of~$\RPC(\omega)$ from
the \emph{mere existence} of such a family.
\end{question}


\begin{question}
Describe the retracts of $\betaN$ and $\Nstar$, as well as their
\emph{absolute} retracts.
\comments
A retract of $\betaN$ is necessarily a closed separable extremally
disconnected subspace.
It is known that a compact separable extremally disconnected can be embedded
as a retract of~$\betaN$.
If~$X$ is such a space then there is a continuous surjection $f:\betaN\to X$
and if $K$~is such that $f\restr K$ is irreducible then $f\restr K$~is
a homeomorphism to~$X$ and $(f\restr K)^{-1}\circ f$ is a retraction
of~$\betaN$ onto~$K$.

Shapiro \cite{MR810825} and
Simon \cite{MR0869226} have shown independently and by quite different means
that not every closed separable subset of~$\betaN$ is a retract.
This gives rise to the notion of an absolute retract of~$\betaN$: a (sub)space
that is a retract irrespective of how it is embedded.

Bella, B\l{}aszczyk and Szyma\'nski proved in~\cite{MR1295157}
that if $X$ is compact, extremally disconnected, without isolated points
and of $\pi$-weight~$\aleph_1$ or less then $X$~is an absolute retract for
extremally disconnected spaces iff $X$~is the absolute of one of the following
three spaces: the Cantor set, the Cantor cube $\vphantom{2}^{\omega_1}2$,
or the sum of these two spaces.
This shows that under~$\CH$ there are very few absolute retracts of~$\betaN$.

We have less information about the retracts of~$\Nstar$, absolute or not.
Of course if a subset of~$\Nstar$ is a retract of~$\betaN$ then it is a retract
of~$\Nstar$ as well.
We do not know whether the converse is true, for separable closed subsets
of course.

We do know that non-trivial zero-sets are not retracts.
Such a set is of the form $Z=\Nstar\setminus\bigcup_{n\in\omega}A_n^*$, where
the $A_n$ are infinite and pairwise disjoint subsets of~$\N$.
We write $C=\bigcup_{n\in\omega}A_n^*$.
Now the closure of~$C$ is a $P$-set in~$\Nstar$, it is the union of~$C$
and the boundary of~$Z$, and if we take one
point~$u_n\in A_n^*$ for each~$n$ then $K=\cl\{x_n:n\in\omega\}$ is a copy
of~$\betaN$ and $K^*=K\setminus\{x_n:n\in\omega\}$ is a $P$-set in the boundary
of~$Z$ and hence in~$Z$.
If we now take assume $r:\Nstar\to Z$ is a retraction then $r\restr K^*$~is
the identity and for all but finitely many~$n$ we must have $r(x_n)\in K^*$.
But this would imply that $K^*$ is separable, a contradiction.

In addition the closure of a non-trivial (not itself closed) cozero-set
may, under~$\CH$ (\cite{MR0248057}), or may not, 
in the $\aleph_2$ Cohen model (\cite{MR0813288}*{Theorem~4.5}), be a retract
of~$\Nstar$. 
\end{question}

\section{Individual Ultrafilters}

\question
Is there a model in which there are no $P$-points and no $Q$-points?
\label{PandQ}
\comments
The Continuum Hypothesis implies that both kinds of points exist.
If $\cee=\aleph_2$ then at least one kind exists; this depends on the
value of~$\dee$.
If $\dee=\cee$ then $P$-points exist, in fact Ketonen showed in~\cite{MR433387}
that then every filter of cardinality less than~$\cee$ can be extended
to a $P$-point.
In the present case, if $\dee<\cee$ then $\dee=\aleph_1$ and then
the result of Coplakova and Vojt\'a\v{s} from~\cite{MR863903} applies
to show that there are $Q$-points; this relies on the fact that the Nov\'ak
number of~$\Nstar$ is at least~$\aleph_2$, see~\cite{MR0600576}.

The current methods for creating models without $P$-points involve iterations
with countable supports and these invariably produce models
where $\cee=\aleph_2$, and hence these will contain $Q$-points.
A recent exception is~\cite{MR3990958}, where models without $P$-points
and arbitrarily large continuum are constructed.
However $\dee=\aleph_1$ in these models, hence these contain $Q$-points
as well.
\endquestion

\begin{question}
Is there a model in which there is a rapid
ultrafilter\pri{rapid ultrafilter|see{ultrafilter, rapid}}%
\pri{ultrafilter!rapid} but in which
there is no $Q$-point\pri{Q-point@$Q$-point}%
\pri{ultrafilter!$Q$-point|see{$Q$-point}}%
\index{ultrafilter!$P$-point|see{$P$-point}}%
\index{ultrafilter!$Q$-point|see{$Q$-point}}%
\pri{ultrafilter!$P$-point|see{$P$-point}}?
\comments
In~\cite{MR4246814} it was shown that the existence of a countable non-discrete
extremally disconnected group implies the existence of rapid ultrafilters.
\end{question}



\begin{question}
What are the possible compactifications of spaces of the form $\N\cup\{p\}$
for $p\in\Nstar$?
\comments
Of course for every $p$ we have $\beta(\N\cup\{p\})=\betaN$.
There are points where this phenomenon persists: Dow and Zhou showed that
is $f:\betaN\to{}^\cee2$ is continuous and onto and $K\subset\Nstar$ is a
closed set such that $f\restr K$ is irreducible and onto then
for \emph{every} point in~$K$ \emph{every} compactification of~$\N\cup\{p\}$
contains a copy of~$\betaN$, see~\cite{MR1676677}.

Other examples of spaces of the form~$\N\cup\{x\}$, where $x$~is the only
non-isolated point, for which every compactification contains~$\betaN$
were constructed by Van~Douwen and Przymusi\'{n}ski in~\cite{MR532957}.

\smallskip
The case of scattered compactifications has received considerable interest.

In~\cite{MR107849} Semadeni asked whether $\N\cup\{p\}$ always has a scattered
compactification.

In~\cite{MR263030} Ryll-Nardzewski and Telgarsky proved that the answer
is yes if $p$~is a $P$-point and the Continuum Hypothesis holds;
the compactification is a version of the compactification~$\gamma\N$
of Franklin-Rajagopalan from~\cite{MR283742},
where $\gamma\N\setminus\N$ is a copy of the
ordinal~$\omega_1+1$ and $p$~corresponds to the point~$\omega_1$.

In~\cite{MR410671} Jayachandran and Rajagopalan constructed a scattered
compactification of $\N\cup\{p\}$, where $p$~is a $P$-point limit
of a sequence of $P$-points.

Solomon, Telgarski, and Malykhin, in \cite{MR448298}, \cite{MR461441},
 and~\cite{MR478101}, respectively, exhibited points~$p$ in~$\Nstar$
such that $\N\cup\{p\}$ has no scattered compactification.

Malykhin's paper and the paper~\cite{MR461442} by Telgarsky contain
investigations of the structure of the (complementary)
sets~$S$ and~$\mathit{NS}$ of points for which $\N\cup\{p\}$ does and
does not have
a scattered compactification respectively.
The set~$\mathit{NS}$ is quite rich: it contains the closures of all of its
countable subsets and it is upward closed in the Rudin-Frol\'ik order.

This richness foreshadowed a later result of Malykhin's
from~\cites{Malykhin:betaNandnotCHa,Malykhin:betaNandnotCHb}:
in the Cohen model it is the case that for \emph{every} point
$p\in\Nstar$ \emph{every} compactification of~$\N\cup\{p\}$
contains a copy of~$\betaN$; in particular $\mathrm{NS}=\Nstar$ in this model.
\end{question}

\begin{question}\label{nonremote}%
Is there $p\in\Q_d^*$ such that
$\calB=\{A\in p:A$~is closed and nowhere dense in~$\Q$ and without
isolated points$\}$ is a base for $p$?
\comments
To eliminate possible confusion: we wrote $p\in\Q_d^*$ to emphasize that
we are asking for an ultrafilter on the countable \emph{set of rationals}
(with the discrete topology),
and $\Q$ in the description of~$\calB$ to emphasize that we want a base for
the ultrafilter that is closely connected to the topological structure of the
\emph{space of rationals}.

One could ask the question in the opposite direction: is there a point~$x$
in $\beta\Q\setminus\Q$ (the \emph{space} of rationals) that,
when considered as an ultrafilter of closed sets has a base consisting
of closed nowhere dense copies of~$\Q$ and that also generates a real
ultrafilter on the set~$\Q$.

A third way of looking at this question is to consider 
$\beta\id:\beta\Q_d\to\beta\Q$, where $\id$ is the identity map and
look for points in~$\beta\Q\setminus\Q$ with one-point preimages.
Such points are easily found in the closure of~$\N$ for example, but we want
a point whose elements are topologically as rich as possible.

These ultrafilters were dubbed `gruff ultrafilters' by Van Douwen.
This question is still open but there are many consistent positive answers:

\begin{itemize}
\item Van Douwen \cite{MR1192307}: from $\MA_{\mathrm{countable}}$,
\item Coplakova and Hart \cite{MR1676672}: from $\bee=\cee$,
\item Ciesielski and Pawlikowski \cite{MR1997781}: from a version of
      the Covering Property Axiom (hence in the Sacks model),
\item Mill\'an \cite{MR2182932}:
       from the same assumption a $Q$-point with this property,
\item Fern\'{a}ndez-Bret\'{o}n and Hru\v{s}\'{a}k \cite{MR3539743}:
      from a parametrized $\diamond$-principle,
      from $\dee=\cee$,
      and in the random real model;
      a correction in~\cite{MR3712981} points out that in the third case
      one needs to add $\aleph_1$~many Cohen reals first
\end{itemize}
\end{question}

\begin{question}
Is there a $p\in\Nstar$ such that whenever $\langle x_n:n\in\omega\rangle$ is
a sequence in $\Q$ there is an $A\in p$ such that $\{\,x_n:n\in A\,\}$
is nowhere dense?
\comments
Such ultrafilters are called \emph{nowhere dense}.
A $P$-point is nowhere dense: it will have a member~$A$ such that
$\{x_n:n\in A\}$ converges to a point or is closed and discrete.
On the other hand, in~\cite{MR1690694} Shelah showed that it is consistent
that there are no nowhere dense ultrafilters.
In~\cite{MR1833478} it is shown that a nowhere dense ultrafilter exists
iff there is a $\sigma$-centered partial order that does \emph{not}
add a Cohen real.

Research into this type of problem was initiated by Baumgartner
in~\cite{MR1335140}:
the general situation involves a set~$S$ and a notion of smallness on~$S$,
usually expressed in terms of ideals.
One then calls an ultrafilter $u$ on~$\N$ small if for every map $f:\N\to S$
there is a member of~$u$ whose image under~$f$ is small.
\end{question}

\begin{question}
Is there an ultrafilter~$u$ such that for every map $f:\N\to\N$ there is
a member~$U$ of~$u$ such that $f[U]$~has density zero?
\comments
This is a special case of the general problem mentioned in the comments above.
We mention it here because it is related to some special cases of
problem~\ref{nwdperm}, which deals with permutations, rather than arbitrary
maps.
\end{question}

\begin{question}
Is there in $\ZFC$ an ultrafilter that is Sacks-indestructible?
\comments
This question is inspired by the many proofs that ultrafilters of small
character may exist.
Sacks forcing preserves selective ultrafilters, $P$-points and many
ultrafilters constructed from these.
Those ultrafilters need not exist of course, so the question becomes
if there are ultrafilters that are preserved by this partial order.
\end{question}

\section{Dynamics, Algebra, and Number Theory}

\begin{question}\label{q.max.orbitcl}%
Is there a point in $\Nstar$ that is not an
element of any maximal orbit closure\pri{orbit closure!maximal}?
\comments
In this problem we consider the integers $\Z$ rather than~$\N$ and the
shift map~$\sigma$, defined by $\sigma(n)=n+1$.
The orbit of~$u\in\Nstar$ is the set $\{\sigma(u):n\in\Z\}$ and its
closure~$C_u$ is the \emph{orbit closure} of~$u$.
\end{question}

\begin{question}
Is there an infinite strictly increasing sequence of
orbit closures\pri{orbit closure}?
\comments
This problem is related to the previous problem: if there is no increasing
sequence of orbit closures then the family of orbit closures is well-founded
under reverse inclusion and every point is in some maximal orbit closure.
A negative answer to this question,
and hence to Question~\ref{q.max.orbitcl},
was given recently by Zelenyuk in~\cite{MR4456336}.
\end{question}

\begin{question}
Is there a $p\in\Nstar$ such that for every pair of commuting
continuous maps $f,g:\Cantor\to\Cantor$ there is an $x\in\Cantor$
such that $\plim f^n(x)=\plim g^n(x)=x$?
\comments
This question is related in two ways to Birkhoff's multiple recurrence
theorem, which states that commuting continuous self-maps of the Cantor set
have common recurrent points.
Using ultrafilters one can state this theorem as: for every pair of commuting
continuous maps $f,g:\Cantor\to\Cantor$ there are~$p\in\Nstar$ and
$x\in\Cantor$ such that $\plim f^n(x)=\plim g^n(x)=x$.

So the first connection to our question is clear: is there one single
ultrafilter that works for all pairs.

The second connection is the question whether the theorem holds
for the Cantor cube~$\Cantorc$?

If it does then the answer to our question is positive.
To see this note first that there are $\cee$~many pairs
of commuting self-maps of~$\Cantor$,
enumerated these as $\{\orpr{f_\alpha}{g_\alpha}:\alpha<\cee\}$.
These determine one pair $\orpr fg$ of commuting self maps of~$\Cantorc$:
write $\Cantorc$ as $\Cantor[\cee\times\omega]$, and let
$f=\prod_{\alpha<\cee}f_\alpha$ and $g=\prod_{\alpha<\cee}g_\alpha$.
The maps~$f$ and~$g$ commute and if $x\in\Cantor[\cee\times\omega]$~is a
common recurrent point then $\plim f^n(x)=\plim g^n(x)=x$ for some~$p\in\Nstar$.
But then also $\plim f_\alpha^n(x_\alpha)=\plim g_\alpha^n(x_\alpha)=x_\alpha$ for
all~$\alpha$.

\end{question}

\begin{question}\label{nwdperm}%
For what nowhere dense sets\pri{nowhere dense set}
$A\subseteq\Nstar$ do we have $\bigcup_{\pi\in S_\N}\pi^*[A]\neq\Nstar$?
\comments
Here $S_\N$ denotes the permutation group of~$\N$.

It is consistent to assume that this happens for all nowhere dense sets.
In~\cite{MR0600576} Balcar, Pelant and Simon studied~$\enn$, the Nov\'ak number
of~$\Nstar$, defined as the smallest number of nowhere dense sets needed to
cover~$\Nstar$.
The inequality~$\cee<\enn$ is consistent and yields the consistency
of ``for all nowhere dense sets'';
it follows from~$\CH$ (because $\enn\ge\aleph_2$),
but is also consistent with other values of~$\cee$.

The inequality $\enn\le\cee$ is also consistent and that case there is not such
an easy way out and it becomes an interesting project to investigate whether
the permutations of individual nowhere sets do, or do not, cover~$\Nstar$
in~$\ZFC$.

Permuting a singleton will not yield a cover, as $|\Nstar|=2^\cee$.

Less obvious is Gryzlov's result from~\cite{MR782711} that the permutations
of the set~$F_d$ from Question~\ref{q.P-sets} do not form a cover.
This was improved by Fla\v{s}kov\'a in~\cite{MR2337416}: the permutations
of the larger set~$F_\Sigma$ do not cover~$\Nstar$ either.

There is another natural nowhere dense subset of~$\Nstar$ the permutations
of which may, or may not, cover~$\Nstar$.
Identify $\N$ with $\N\times\N$ and for~$k\in\N$ and $f:\N\to\N$ write
$U(f,k)=\{\orpr mn:m\ge k$ and~$n\ge f(m)\}$.
The set $B=\bigcap_{f,k}U(f,k)^*$ is nowhere dense and it is well known
then $\bigcup_{\pi\in S_\N}\pi[B]$ consists of all non $P$-points of~$\Nstar$.
Hence the permutations of~$B$ cover~$\Nstar$ iff there are no $P$-points.
\end{question}

\begin{question}
For what nowhere dense sets\pri{nowhere dense set}
$A\subseteq\Nstar$ do we have
$\bigcup\{h[A]:h\in\Aut\}\neq\Nstar$?
\comments
This question is more difficult than the previous one.

For example, singleton sets still do not provide covers in~$\ZFC$, but the
easy counting argument is replaced by the non-trivial fact that $\Nstar$~is not
homogeneous.

We have no information about the sets $F_d$ and $F_\Sigma$ in this context,
except for the general fact that under~$\CH$ the space~$\Nstar$ cannot be
covered by nowhere dense $P$-sets, see~\cite{MR548097}.
Also, in~\cite{MR620204} it was shown that it is consistent that
$\Nstar$~can be covered by nowhere dense $P$-sets,
and the principle~$\NCF$ (Near Coherence of Filters) implies that
$\Nstar$~is even the union of a \emph{chain} of nowhere dense $P$-sets,
see~\cite{MR1261700},
but the sets in these covers are unrelated to the sets~$F_d$ and~$F_\Sigma$.
It is also unclear whether any one of the individual sets in these families
will produce a cover when moved around by the members of~$\Aut$.

The answer for the set~$B$ remains the same because
the union $\bigcup\{h[B]:h\in\Aut\}$ consists of all non-$P$-points.
\end{question}

\section{Other}

\begin{question}\label{q.Katowice}%
Are $\omega_0^*$ and $\omega_1^*$ ever homeomorphic%
\pri{homeomorphism!of $\omega^*$ and $\omega_1^*$}?
\comments
This is known as the Katowice Problem, or rather the last remaining case of
this problem.
It was posed in full by Marian Turza\'nski, when he was in Katowice
(hence the name of the problem).
The general question is:
if $\kappa$ and~$\lambda$ are infinite cardinals, endowed with the discrete
topology, and the remainders $\kappa^*$ and~$\lambda^*$ are homeomorphic
must the cardinals~$\kappa$ and~$\lambda$ be equal?

Since the weight of $\kappa^*$ is equal to~$2^\kappa$ it is immediate that
the Generalized Continuum Hypothesis implies a yes answer.
In joint work Balcar and Frankiewicz established that the answer is actually
positive without any additional assumptions,
\emph{except possibly for the first two infinite cardinals}.
More precisely, see~\cites{MR0461444,MR511955}:
If $\orpr\kappa\lambda\neq\orpr{\aleph_0}{\aleph_1}$
and $\kappa<\lambda$ then the remainders $\kappa^*$ and~$\lambda^*$
are not homeomorphic.

The paper~\cite{MR3563083} contains a list of the current known of consequences
of~$\omega_0^*$ and~$\omega_1^*$ being homeomorphic; all but one of these
can be made to hold in a single model of~$\ZFC$.

By Stone-duality the Katowice problem can be formulated algebraically:
are the quotient (Boolean) algebras~$\pow(\omega_0)/\fin$
and~$\pow(\omega_1)/\fin$ ever isomorphic?
In this form the question even makes sense in~$\ZF$: in models without
non-trivial ultrafilters the spaces~$\omega_0^*$ and~$\omega_1^*$ are empty
(and so trivially homeomorphic) but the structures of the algebras may still
differ.
\end{question}

\begin{question}
Is there consistently an uncountable cardinal $\kappa$ such that
$\omega^*$ and $U(\kappa)$ are homeomorphic?
\comments
This problem is part of the uniform version of the Katowice problem,
Question~\ref{q.Katowice}.
The full question asks whether for distinct infinite cardinals $\kappa$
and~$\lambda$ spaces $U(\kappa)$ and $U(\lambda)$ of \emph{uniform}
ultrafilters can be homeomorphic, or algebraically whether the quotient
algebras $\pow(\kappa)/[\kappa]^{<\kappa}$
and~$\pow(\lambda)/[\lambda]^{<\lambda}$ can be isomorphic.
This is Question~47 in~\cite{MR588216}, where we also find the information
that, in general, the algebra $\pow(\kappa)/[\kappa]^{<\kappa}$ has
cardinality~$2^\kappa$ and is $\mu$-complete for~$\mu<\cf\kappa$ but
not $\cf\kappa$-complete.
Therefore we can concentrate on cases where $2^\kappa=2^\lambda$ and
$\cf\kappa=\cf\lambda$.

In~\cite{MR1103989} Van Douwen investigated the statements~$S_n$:
$$
\hbox{if $\kappa\neq\aleph_n$ then $\pow(\kappa)/[\kappa]^{<\kappa}$ and
$\pow(\omega_n)/[\omega_n]^{<\aleph_n}$ are not isomorphic}.
$$
Thus, our question is whether it is consistent that $S_0$ is false.
Van Douwen showed that there is at most one~$n$ for which $S_n$~is false,
but the proof offers no information on the location of that~$n$ (if any)
as it simply establishes the implication
``if $m<n$ and $S_m$~is false then $S_n$~holds''.
\end{question}

\begin{question}
What is the structure of the sequences
$\langle n((\Nstar)^n): n\in\N\rangle$
and $\langle wn((\Nstar)^n): n\in\N\rangle$?
\comments
Here $n$ and $wn$ denote the Nova\'k and weak Nov\'ak numbers, defined
as the minimum cardinality of a family of nowhere dense sets that covers
the space, or has a dense union, respectively.

It is clear that if $N$~is nowhere dense in a space~$X$ then $N\times Y$
is nowhere dense in the product~$X\times Y$.
This shows that, in general, $n(X\times Y)\le\min\{n(X),n(Y)\}$ and likewise
for~$wn$.
It follows that both sequences in our question are non-increasing and hence
must become constant eventually.

One could ask when they do become constant.
For $wn$ this is undetermined: in~\cite{MR1641157} Shelah and Spinas showed
that for every~$n$ there is a model in
which $wn((\Nstar)^n)> wn((\Nstar)^{n+1})$.
In particular $wn(\Nstar)>wn(\Nstar\times\Nstar)$ is possible,
in~\cite{MR1751223} the latter inequality was shown to hold in the Mathias
model.

For the Nov\'ak numbers of the finite powers nothing is known as yet.
\end{question}

\begin{question}
What is the status of the statement
that all Parovichenko spaces are co-absolute (with~$\Nstar$)?
\comments
This question is related to Parovichenko's theorem from~\cite{MR0150732},
which states that under~$\CH$ all Parovichenko spaces are homeomorphic
to~$\Nstar$.
Of course Parovichenko spaces were named after this theorem was proved:
they are compact, zero-dimensional $F$-spaces of weight~$\cee$ without isolated
points in which every non-empty $G_\delta$-set has non-empty interior.
For the nonce we say that a space is of \emph{Parovichenko type} if it
satisfies the conditions above, except for possibly the weight restriction.

In \cite{MR612009} Broverman and Weiss proved that under~$\CH$ all spaces
of Parovichenko~type of $\pi$-weight~$\cee$ are co-absolute (with~$\Nstar$).
They also established that if~$\CH$ fails and $\cee=2^{<\cee}$ then there
is a Parovichenko space that is not co-absolute with~$\Nstar$.
They also proved that $\omega_0^*$ and~$\omega_1^*$ are co-absolute or,
in algebraic terms that the Boolean algebras~$\pow(\omega_0)/\fin$
and~$\pow(\omega_1)/\fin$ have isomorphic completions, which shows that
completions do not have a direct effect on Question~\ref{q.Katowice}.

In~\cite{MR648079} Williams also established the $\pi$-weight result and
showed that $\Nstar$ is co-absolute with a linearly ordered space.

In~\cite{MR676966} Van Mill and Williams improved the negative result
of Broverman and Weiss: if our statement holds then not only do we
have $\cee<2^{<\cee}$, but even $\cee<2^{\aleph_1}$.

In~\cite{MR759135} Dow proved that the equality $\cf\cee=\aleph_1$
already implies that all Parovichenko spaces are co-absolute.

The definition of the absolute as the Stone space of the Boolean algebra
of regular open sets makes sense for any compact space, so one may also seek
co-absolutes of~$\Nstar$ among spaces that are not zero-dimensional.
In~\cite{MR234422} Comfort and Negrepontis showed that under~$\CH$ if
$X$~is locally compact and $\sigma$-compact, but not compact, and if
$\bigl|C(X)\bigr|=\cee$ then the set of $P$-points in~$X^*$ is homeomorphic
to the $G_\delta$-modification of the ordered space~$2^{\omega_1}$;
Parovichenko had already established this fact for~$\Nstar$ in~\cite{MR0150732}.
This implies that for such spaces the remainders share a homeomorphic
dense subspace and hence that all such remainders are co-absolute
with~$\Nstar$, still under~$\CH$ of course.
So, for example, under~$\CH$ the spaces~$\Nstar$ and $\Hstar$ are co-absolute.

In~\cite{MR1619290} Dow showed that in the Mathias model $\Nstar$ and~$\Hstar$
are not co-absolute.
\end{question}

\begin{question}
Let $X$ be a compact space\pri{compact space} that can be mapped
onto~$\Nstar$.
Is $X$ non-homogeneous\pri{compact space!non-homogeneous}?
\comments
Since $\Nstar$ maps onto~$\betaN$, a space as in the question will also
map onto~$\betaN$.
If the weight of~$X$ is at most~$\cee$ then Theorem~4.1\,(c) of~\cite{MR644652}
applies and we find that $X$~is indeed non-homogeneous.
\end{question}

\begin{question}
Is it consistent that every compact space\pri{compact space}
contains either a convergent sequence\pri{sequence!convergent}
or a copy of~$\betaN$\pri{copy!of~$\beta\omega$}?
\label{q.efimov}
\comments
Efimov asked in~\cite{MR0253290} whether every compact space contains
either a convergent sequence or a copy of~$\betaN$ and a counterexample
is now called a \emph{Efimov space}.
In~\cite{hart:efimov} one finds a survey of the status of the problem
in~2007; it lists various consistent Efimov spaces, which explains
why the present formulation asks for a consistency result.
We mention here some of the additional results that have been obtained
in the meantime.

To begin there is a positive answer in~\cite{MR3164725} to Question~1
from~\cite{hart:efimov}:
Martin's Axiom, or even the equality $\bee=\cee$, implies that there is
a Efimov space.

In addition there has been progress on two related questions due to Juh\'asz
and Hu\v{s}ek.
The latter asked whether every compact Hausdorff space contains either
a convergent $\omega$-sequence or a convergent $\omega_1$-sequence; Juh\'asz'
question is stronger: must a compact Hausdorff space that does not contain
a convergent $\omega_1$-sequence be first-countable?
A counterexample to Hu\v{s}ek's question would be a Efimov space because
$\betaN$~contains a convergent $\omega_1$-sequence.
In~\cite{MR1707489} one finds a result that provides many models in which
Juh\'asz' question, and hence that of Hu\v{s}ek's, have a positive answer.
One of these models satisfies $\bee=\cee$, hence Efimov's question is strictly
stronger than that of Hu\v{s}ek's.
\end{question}

\begin{question} Is there a locally connected continuum such that every proper
          subcontinuum contains a copy of~$\betaN$?
\comments
There are various continua that have the property that every proper
subcontinuum contains a copy of~$\betaN$: the remainders
$\beta\R^n\setminus\R^n$ all have this property for example.
The reason is that they are $F$-spaces, hence the closure of every countable
relatively discrete subset is a copy of~$\betaN$.
However, these remainders are not locally connected; indeed if a space~$X$
is not pseudocompact then one can use an unbounded continuous function
to exhibit points in~$X^*$ at which neither $\beta X$ nor $X^*$~is locally
connected, see~\cite{MR96195}.\par
In~\cite{MR1934264} we find a construction, from~$\CH$, of a locally connected
continuum without non-trivial convergent sequences.
This construction, an inverse limit in which all potential convergent
sequences are destroyed, can be modified with some extra bookkeeping to yield
a locally connected continuum in which every infinite subset contains
a countable discrete subset whose closure is homeomorphic to~$\betaN$,
still under~$\CH$ of course.\par
This leaves the question for a $\ZFC$-example open but also suggest some further
variations.
The example has the property that \emph{some} countable relatively discrete
subsets have $\betaN$~as their closures.
One can ask whether one can ensure this for \emph{all} countable relatively
discrete subsets, or whether one can even make all countable subsets
$C^*$-embedded.
The reason for this is that a compact $F$-space cannot be locally connected,
hence we would like to know how close to an $F$-space a locally connected continuum
can be.\par
We would also like to know whether there is a natural example that answers
our question; natural in the sense that one can simply write it down, as in
``$\betaN$~is a compact space without convergent sequences'' and
``$\Hstar$~is a continuum in which every proper subcontinuum contains
  a copy of~$\betaN$''.
\end{question}

\begin{question}
Is there an extremally disconnected normal locally compact space%
\pri{extremally disconnected!normal!locally compact}%
\pri{locally compact space!normal!extremely disconnected}
          that is not paracompact\pri{compact space!not paracompact}?
\comments
The ordinal space $\omega_1$ is locally compact and normal, but not paracompact.
There are, however, various additional assumptions that when added
to local compactness and normality will ensure paracompactness.
Extremal disconnectedness may or may not be such an assumption:
Kunen and Parsons showed in~\cite{MR540504} that if $\kappa$~is weakly compact
then $\beta\kappa\setminus U(\kappa)$~is normal and locally compact but not
paracompact.
As weak compactness is a large cardinal property the answer to this question
can go many ways: a consistent counterexample, a real counterexample, or
even an equiconsistency result involving a large cardinal.

The weaker property of basic disconnectedness does not work, as shown
by Van Douwen's example in~\cite{MR546947}.
In this paper Van Douwen attributes the present question to Grant Woods.
\end{question}

\begin{question}
Is every compact hereditarily paracompact space of weight at most~$\cee$ a
continuous image of~$\Nstar$?
Is every hereditarily c.c.c.\ compact space a continuous image of~$\Nstar$?
\comments
These questions are part of the general problem of identifying
the continuous images of~$\Nstar$.
Przymusi\'{n}ski proved in~\cite{MR671232} that all perfectly normal
compact spaces are continuous images of~$\Nstar$.
One can therefore look for weakenings of perfect normality that still
make the space an image of~$\Nstar$.
The present two properties are such weakenings and they have
not been ruled out yet.

Another weakening, first-countability, was ruled out by Bell
in~\cite{MR1058795}: the $\aleph_2$-Cohen model contains a first-countable
compact space that is not a continuous image of~$\Nstar$; this space is
also hereditarily metacompact.
In the same paper Bell showed that the compact ordered space~$2^{\omega_1}$
(with the lexicographic order) is an image of~$\Nstar$.
Theorems~15 and~17 in Chapter~1 of~\cite{MR0220252} imply that every
compact ordered space that is first-countable is a continuous image of the
latter space, hence also of~$\Nstar$.

In connection with the latter result we note that it is consistent with
the negation of~$\CH$ that all linear orders of cardinality~$\cee$ are
embeddable into the Boolean algebra~$\powNfin$, see~\cite{MR567675}.
By a combination of the Stone and Wallman dualities this implies that
it is consistent with $\neg\CH$ that every compact ordered space of
weight~$\cee$ is a continuous image of~$\Nstar$.

This was later generalized in~\cite{MR1095691} to the consistency of
Martin's Axiom for $\sigma$-linked partial orders, the negation of~$\CH$,
and the statement that all compact spaces of weight~$\cee$ are continuous
images of~$\Nstar$.

In both cases the proof constructs an embedding of a universal linear order
or a universal Boolean algebra of cardinality~$\cee$ into~$\powNfin$.
This raises the question whether there is a universal compact space
of weight~$\cee$; one that maps onto all such spaces.
The answer is negative, see~\cite{MR1707489}*{Section~6}.
\end{question}

\begin{question}
Is every compact space of weight at most~$\aleph_1$ a $1$-soft remainder
of~$\omega$?
\comments
A compactification~$\gamma\N$ of~$\N$ is $1$-soft if for every subset~$A$
of~$\N$ with $\cl A\cap\cl(\N\setminus A)\neq \emptyset$ there is an
autohomeomorphism~$h$ of~$\gamma\N$ that is the identity
on~$\gamma\N\setminus\N$ and is such that $\{n\in A:h(n)\notin A\}$ is
infinite.

See Question~351527 on MathOverFlow, \cite{flow351527},
and also the papers \cite{MR4142223} and~\cite{MR4266612} for related
information.
\end{question}

\begin{question}
Is there a universal compact space of weight~$\aleph_1$?
\comments
We mean universal in the mapping-onto sense;
the dual question has the well-known answer~$[0,1]^{\omega_1}$ and Parovichenko's
theorem suggests that the answer might be positive.
The answer is negative in the $\aleph_2$-Cohen model but a good reference
is hard to find.
There are references to the result,
\cites{MR1676311,MR739914,MR1057268,MR1056364},
but no concrete proof.

However, the argument in~\cite{MR1707489}*{Section~6} can readily be adapted
to provide an accessible proof.
We apply Stone duality and show that in the model there is no Boolean algebra
of cardinality~$\aleph_1$ in which every Boolean algebra of that cardinality
can be embedded.
Let $\Fn(\omega_2\times\omega_0,2)$ denote the Cohen partial order and
let $G$ be a generic filter.

The main steps are: we can assume that the Boolean algebra is determined
by a partial order~$\prec$ on a subset of~$\omega_1$.
By the ccc of $\Fn(\omega_2\times\omega_0,2)$ the order~$\prec$ is a member
of~$V[G\restr\alpha]$ for some~$\alpha<\omega_2$.
Take the next $\aleph_1$ many Cohen reals
$\langle c_\beta:\beta<\omega_1\rangle$,
defined by $c_\beta(n)=\bigcup G(\alpha+\beta,n)$.
The union, $T$, of the binary tree $2^{<\omega}$ and the set
$\{c_\beta:\beta<\omega_1\}$ is a partially ordered set which, when turned
upside-down generates a Boolean algebra~$B$.
Assume $\varphi: T\to\omega_1$ is the restriction of an embedding of~$B$
into $\orpr{\omega_1}{\prec}$.
There is a countable subset~$C$ of~$\omega_2$ such that the restriction
of~$\varphi$ to~$2^{<\omega}$ belongs to~$V[G\restr(\alpha\cup C)]$.
Now take $\beta\in\omega_1$ such that $\alpha+\beta\notin C$.
Then $c_\beta$ does not belong to~$V[G\restr(\alpha\cup C)]$, yet
it can be defined from the elements $\gamma=\varphi(c_\beta)$
and $\varphi\restr2^{<\omega}$ by the
formula $\bigcup\{s:\gamma\prec\varphi(s)\}$.
\end{question}

\begin{question}
Investigate ultrafilters as topological spaces.
\comments
This is a very general question, so let us discuss some specific ones that
may be investigated.
An ultrafilter can be viewed as a subspace of the Cantor set~$\Cantor$,
if one identifies a subset of~$\omega$ with its characteristic function.

Of course this makes ultrafilters separable metric spaces, and hence
relatively well-behaved.
But not too well-behaved: free ultrafilters are non-measurable and do not
have the property of Baire.

To begin one can repeat many of the investigations into the Rudin-Keisler
order using more general kinds of maps.
We know $p\leRK q$ means that there is a map $\varphi:\N\to\N$ such
that $\beta\varphi(q)=p$.
The map~$\varphi$ determines a continuous map from~$\Cantor$ to itself, so
the following definition suggests itself at once: say $p\le_c q$ if there
is a continuous map $f:\Cantor\to\Cantor$ such that $f[q]=p$.

One can ask whether $p\le_cq$ and $q\le_cp$ together imply that $p\equiv_cq$,
which means that there is a homeomorphism of~$\Cantor$ that maps $p$ to~$q$.
The structure of the partial order~$\le_c$, minimal elements, incomparable
elements, etc., would warrant investigation as well.

There is no reason to stop there of course: one can ask the same questions
about Borel maps of any specific order, or of maps of arbitrary Baire classes.

One need not work with maps on~$\Cantor$, though that may make life easier,
one can investigate what it means for two ultrafilters to be homeomorphic,
or what it means that one is a continuous image of the other.
The methods of~\cite{MR1277880} may be of use in determining the possible
sizes of sets of ultrafilters that are incomparable in this sense.

We note that an ultrafilter can be homeomorphic to at most $\cee$ many other
ultrafilters: if $f:p\to q$ is a homeomorphism then Lavrentieff's theorem
implies that $f$ can be extended to a homeomorphism of $G_\delta$-subsets
of~$\Cantor$, and the number of such homeomorphisms is equal to~$\cee$.

The paper~\cite{MR2879361} contains many results on the topology of
ultrafilters.
\end{question}

\begin{question}
Is it consistent that all free ultrafilters have the same Tukey type?
\comments
Isbell~\cite{MR201316} raised the question of the number of Tukey types
of ultrafilters on~$\N$ and gave the obvious bounds $2$ (trivial or not)
and $2^\cee$.
Tukey types of free ultrafilters were investigated
by Dobrinen and Todor\v{c}evi\'c in~\cite{MR3069290}
who gave a combinatorial characterization of ultrafilters that are
Tukey-equivalent to the partial order of finite subsets of~$\cee$:
the ultrafilter~$\calU$ should contain a subfamily~$\calX$ of
cardinality~$\cee$ such that for every infinite subfamily~$\calY$
of~$\calX$ the intersection~$\bigcap\calY$ does not belong to~$\calU$.

Such ultrafilters exist see~\cite{MR201316}*{Theorem~5.4}; they
are the ultrafilters of character~$\cee$ constructed from an independent
family of cardinality~$\cee$, see also~\cite{MR0001454}.

In~\cite{MR3990958}*{Announcement~9} Chodounsk\'y and Guzm\'an announce
a result that comes close to the statement that all free ultrafilters
have this property.

\textbf{Added in proof}:
in~\cite{arxiv:2410.08699} Cancino-Manr\'iquez and Zapletal construct
models where all free ultrafilters are Tukey equivalent
to the partial order of finite subsets of~$\cee$.
\end{question}

\begin{question}
Is the space of minimal prime ideals of $C(\Nstar)$ \emph{not}
basically disconnected?
\comments
For a commutative ring~$R$ we let $mR$ denote the set of minimal prime ideals
endowed with the hull-kernel topology.
In~\cites{MR0144921,MR194880} Henriksen and Jerison studied this space and
asked whether $mC(\Nstar)$~is basically disconnected.

In the papers~\cite{MR958091} and~\cite{MR1057626} various conditions were
found that imply $mC(\Nstar)$~is \emph{not} basically disconnected.
For example, $\MA$~implies that $mC(\Nstar)$ is not even an $F$-space
(\cite{MR958091}).
In~\cite{MR1057626} it was shown that the equality $\cf[\dee]^{\aleph_0}=\dee$
suffices to show that $mC(\Nstar)$~is not basically disconnected.
Failure of this equality entails the existence of inner models with measurable
cardinals.
The actual consequence, called~\textbf{Mel}, of this equality that was used
in the proof identifies $\N$ with~$\Q$ and asks for a $P$-filter~$\calF$
on~$\Q$, and two countable disjoint dense subsets $A$ and $B$
of $\R\setminus\Q$ such that the closure in~$\R$ of every member
of~$\calF$ meets both $A$ and~$B$.

Thus, to show that $mC(\Nstar)$ is not basically disconnected it suffices
to show that \textbf{Mel} holds, or the following stronger, but possibly more
manageable, statement: the ideal of nowhere dense subsets of $\Q$ can be
extended to a $P$-ideal.
\end{question}

\begin{question}
Is there a c.c.c.\ forcing extension%
\pri{forcing extension!\ccc}\pri{ccc\ccc{}!forcing extension}
of~$L$ in which there are no $P$-points\pri{P-point@$P$-point}?
\comments
The consistency of the nonexistence of $P$-points was proven by Shelah,
see~\cite{MR728877} and also~\cite{MR1623206}*{VI\,\S4}.

After this there have been various attempts to (dis)prove the existence
of $P$-points in various standard models.
Quite often the outcome was that ground model $P$-points remained ultrafilters
and $P$-points in the extension.

A notable exception is the Silver model: in~\cite{MR3990958} we find a proof
that iterating Silver forcing $\omega_2$ times with countable supports produces
a model without $P$-points; the same holds for the countable support product
of arbitrarily many copies of the partial order.
This establishes the consistency of the nonexistence of $P$-points
with arbitrarily large values of~$\cee$.

A question that is still open is whether $P$-points exist in the random real
model.
If not then this would answer the present question positively.
If there are $P$-points in this model then our question gains interest
as it is as yet unknown whether c.c.c.\ forcing can be used to kill
$P$-points.
\end{question}

\begin{question}
What is the relationship between ultrafilters of small character
(less than~$\cee$) and $P$-points?
\comments
One of the first ultrafilters of small character can be found
in~\cite{MR597342}*{Exercise~VII.A10}; it is a simple $P_{\aleph_1}$-point
constructed by iterated forcing over a model of~$\neg\CH$.
There are many more examples of ultrafilters of small character but their
constructions seem to involve $P$-points in some form or another.
A common method is to start with a model of~$\CH$ and enlarge the continuum
while preserving some ultrafilters; these will then have character~$\aleph_1$,
which is smaller than~$\cee$.
Almost always these `indestructible' ultrafilters are $P$-points (or stronger)
and remain $P$-points in the extension.
There are a few exceptions, see~\cite{MR987317} for instance, but there the
ultrafilters are built using $P$-points and these are preserved as well.
\end{question}

\begin{question}
We let $\Spchi$ denote the set of characters of ultrafilters on~$\N$,
the \emph{character spectrum} of~$\N$.
The general question is what one can say about this set.
\comments
We know that $\cee\in\Spchi$, and that $\Spchi=\{\cee\}$ is possible.

In \cite{MR2365799} Shelah showed the consistency of there being
three cardinals $\kappa$, $\lambda$, and~$\mu$ such that
$\kappa<\lambda<\mu$, and 
$\kappa,\mu\in\Spchi$ and $\lambda\notin\Spchi$.
The construction uses a c.c.c.~forcing over a ground model in which
the three cardinals are regular, $\lambda$~is measurable, and there is another
measurable cardinal below~$\kappa$.
In~\cite{MR2847327} he extended this result by showing how to build,
given two disjoint sets $\Theta_1$ and~$\Theta_2$ of regular cardinals,
a cardinal-preserving partial order that forces $\Theta_1$ to be a subset
of~$\Spchi$ and $\Theta_2$ to be disjoint from it; the construction requires
$\Theta_2$ to consist of measurable cardinals.
The same paper also contains models in which $\{n:\aleph_n\in\Spchi\}$ can be
any subset of~$\N$, starting from infinitely many compact cardinals.
This answers a question from~\cite{MR1686797}, namely whether if there
are ultrafilters of character~$\aleph_1$ and~$\aleph_3$ there must be one
of character~$\aleph_2$, but at the cost of large cardinals.

This leaves open the question whether the conjunction of
$\aleph_1,\aleph_3\in\Spchi$ and $\aleph_2\notin\Spchi$ can be proven
consistent from the consistency of just~$\ZFC$.
To be very specific we ask whether there is an ultrafilter of
character~$\aleph_2$ in the model(s) of~\cite{MR597342}*{Exercise~VII.A10},
where one starts with a model of~$\cee=\aleph_3$, and in the side-by-side
Sacks model where $\cee=\aleph_3$.
\end{question}


\begin{question}
Is there consistently a point in $\Nstar$ whose
$\pi$-character\pri{pi-character@$\pi$-character} has
countable cofinality?
\comments
The paper~\cite{MR1686797} contains a wealth of material on $\pi$-characters
of ultrafilters, including a model with an ultrafilter of
$\pi$-character~$\aleph_\omega$.

Unlike the results on the character spectrum the results on the $\pi$-character
spectrum do not require large cardinals.
\end{question}

\begin{question}
Is it consistent that $t(p,\Nstar)<\chi(p)$\pri{tpomega@$t(p,\Nstar)$}%
\pri{chip@$\chi(p)$} for some $p\in\Nstar$?
\comments
There are plenty of compact spaces with points where the tightness is smaller
than the character; the one-point compactification of the any uncountable
discrete space will do: the tightness at the point at infinity is countable,
the character of the point is not.

Let us remark that no point of $\Nstar$ has countable tightness: certainly
at $P$-points the tightness is uncountable; if $p$~is not a $P$-point
then it lies on the boundary of a zero-set~$C$ and in the closure of its
interior, but the closure of every countable subset of that interior
is a subset of that interior.
This implies that $t(p,\Nstar)=\chi(p)=\cee$ if~$\CH$ holds, hence the
question for a consistency result.

As an aside we mention that there are consistent examples
of regular extremally disconnected spaces of countable tightness:
in~\cite{MR0458392} and~\cite{MR461464} one finds constructions of
extremally disconnected $S$-spaces.
The constructions use $\clubsuit$ and that some extra assumption is necessary
is shown in~\cite{MR615971}: there are no extremally disconnected $S$-spaces
if $\MAnotCH$~holds.
Both \cite{MR461464} and~\cite{MR615971} contain constructions of extremally
disconnected $S$-spaces in~$\betaN$.
\end{question}

\begin{question}
If $C(\omega+1,\C)$\pri{Comega+1@$C(\omega+1)$} admits an
incomplete\pri{norm!incomplete} norm then does
$C(\betaN,\C)$\pri{Cbetaomega@$C(\beta\omega)$} admit one too?
\label{analyse.1}
\comments
This question is related to a conjecture\slash question of Kaplansky's about
algebra norms on the spaces~$C(X,\C)$, with $X$ compact.
The question is whether every algebra norm is equivalent to the
sup-norm~$\|{\cdot}\|_\infty$.
The answer is positive if the norm is complete, hence the question became
whether every algebra norm on~$C(X,\C)$ is complete.

The book~\cite{MR942216} surveys the solution to this problem:
under~$\CH$ every~$C(X,\C)$ carries an incomplete algebra norm
(Dales and Esterl\'e) and it is consistent
that every algebra norm on every~$C(X,\C)$ is complete (Solovay and Woodin).

The present question comes from the results that if $C(\betaN,\C)$~admits
an incomplete norm then so does every~$C(X,\C)$, and
if \emph{some}~$C(X,\C)$ carries an incomplete norm then so
does~$C({\omega+1},\C)$.
In short it asks whether all compact spaces are equivalent for Kaplansky's
conjecture.

The question can be translated into terms of individual ultrafilters and this
leads to some interesting subquestions.
A seminorm on an algebra is a function that satisfies all conditions of an
algebra norm except for the condition that non-zero elements should have
non-zero norm.
An algebra is semi-normable if it carries a non-trivial seminorm.

For a point $p$ of $\betaN$ we let $A_p$ denote the quotient algebra
$M_p/I_p$, where $M_p=\{f\in C(\betaN,\C):f(p)=0\}$, and
$I_p=\{f\in C(\betaN,\C):(\exists P\in p)(f\restr P=0)\}$.
We also let $c_0$ be the subalgebra of~$C(\betaN,\C)$ of functions
that vanish on~$\Nstar$ and we let $c_0/p$ denote the quotient algebra
$c_0/(c_0\cap I_p)$.

Theorem~2.21 in~\cite{MR942216} shows why we should be interested
in these algebras:
The algebra $C(\betaN,C)$ admits an incomplete norm iff for some~$p$
the algebra~$A_p$ is seminormable, and $C({\omega+1},\C)$ admits an incomplete
norm iff for some~$q$ the algebra~$c_0/q$ is seminormable.

We see that \emph{if} there is a $p$ such that $A_p$ is seminormable \emph{then}
there is a~$q$ such that $c_0/q$~is seminormable.
The present question ask whether this implication can be reversed.

Further questions regarding these algebras suggest themselves:
is it the case that the seminormability of~$A_p$ implies that of~$c_0/p$?
In other words can we get $q=p$ in the previous paragraph?

Also, what is the answer to the stronger version of our question:
if $c_0/p$~is seminormable is $A_p$ seminormable too?

We recommend \cite{MR942216}*{Chapters~1, 2 and~3} for more detailed
information on this question.
\end{question}

\begin{question}
($\MA+\lnot\CH$\pri{MACH@$\MA+\lnot\CH$}) Are there $\Gap$ and $p$
($P$-point, selective) such that $p\subseteq I_\Gap^+$?
\comments
Here $\Gap$ denotes a Hausdorff-gap in~${}^\omega\omega$ and
$I_\Gap$ is the ideal of sets over which $\Gap$~is filled.

S. Kamo \cite{MR1138202} proved that if $V$~is obtained from a model of~$\CH$
by adding Cohen reals then in~$V$ an ideal is a gap-ideal iff it is
$\le\aleph_1$-generated.
Also, $\CH$~implies that any nontrivial ideal is a gap-ideal.

The commentary in~\cite{MR2023411} mentions a further preprint by Kamo,
\cite{kamo1993}, where it is shown that, under $\MA+\lnot\CH$,
for every Hausdorff gap~$\Gap$ there are both selective ultrafilters and
non-$P$-points consisting of positive sets (with respect to the
gap-ideal~$I_\Gap$).
Also under $\MA + \lnot\CH$ there is a selective non-$P_{\aleph_2}$-point
that meets every gap-ideal.

Unfortunately we were unable to locate this preprint and verify these
statements.
\end{question}

\section{Orders}

\begin{question}
Is there for every $p\in\Nstar$ a $q\in\Nstar$
such that $p$ and $q$ are $\leRK$-incomparable%
\pri{ultrafilters!$\leRK$-incomparable}?
\comments
This question has a long history; it is as old as the Rudin-Keisler order
itself.
In~\cite{MR314619} Kunen constructed two points that are $\leRK$-incomparable.
In~\cite{MR540490} Shelah and Rudin proved that there is a set of $2^\cee$
incomparable points.
In~\cite{MR0863940} Simon proved that these points may be taken to
be $\aleph_1$-OK.
In~\cite{MR1277880} Dow showed that that there are many more situations
where such sets may be constructed.

However, none of these results shed light on the present question.
Some partial results are available: in~\cite{MR931732} Hindman
proved: if $p$ is such that $\chi(r)=\cee$ whenever $r\leRK p$ then there
is a point that is incomparable with~$p$, so the answer to the present question
is positive if all ultrafilters have character~$\cee$.
Furthermore if $\cee$~is singular and $\chi(p)=\cee$ then again there is a
point that is incomparable with~$p$.
The latter result was extended by Butkovi\v{c}ov\'a in~\cite{MR1045131}:
if $\kappa<\cee$ is such that $\cee<2^\kappa$ then for every ultrafilter
of character~$\cee$ there are $2^\kappa$~many ultrafilters incomparable
with it.
Note that these results all impose conditions on individual ultrafilters
in order to find an incomparable point; only the condition
``all ultrafilters have character~$\cee$'' answers this question directly.

In~\cite{MR1623206}*{XVIII~\S4} Shelah proved that it is consistent that
up to permutation there is one $P$-point.
\end{question}

We recall the definition of the Rudin-Frol\'ik order: we say $p\leRF q$ if
there is an embedding $f:\betaN\to\betaN$ such that $f(p)=q$.
This is a preorder that induces a partial order on the types of ultrafilters.
To see this note that $p\leRF q$ implies $p\leRK q$: given~$f$ take a partition
$\{A_n:n\in\N\}$ of~$N$ such that~$A_n\in f(n)$ for all~$n$.
The map $g=\bigcup_n(A_n\times\{n\})$ satisfies $p=g(q)$ and shows $p\leRK q$.

As usual $p\lRF q$ will mean $p\leRF q$ plus not-$q\leRF p$, and this is
readily seen to be equivalent to there being an embedding $f:\betaN\to\Nstar$
such that~$f(p)=q$.

The Rudin-Frol\'ik order is tree-like: if $p,q\leRF r$ then $p\leRF q$
or $q\leRF p$.
And due to the relation with~$\leRK$ we see at once that $\{p:p\leRK q\}$
always has cardinality at most~$\cee$.

In many papers on the Rudin-Frol\'ik order Frol\'ik's original notation
is employed where one writes $X=f[\N]$, and $q=\Sigma(X,p)$ as well
as $p=\Omega(X,q)$.

\begin{question}
For what cardinals $\kappa$ is there a strictly decreasing chain
of copies of~$\betaN$ in $\Nstar$ with a one-point intersection?
\label{q.one-point}
\comments
This question is related to decreasing chains in~$\leRF$.
A decreasing sequence of copies of~$\betaN$ determines and is determined
by a sequence $\langle X_\alpha:\alpha\in\delta\rangle$ of countable
discrete subsets of~$\Nstar$ with the property
that $X_\alpha\subseteq\cl{X_\beta}\setminus X_\beta$ whenever $\beta\in\alpha$.
Take a point~$p$ in the intersection of the sequence;
then $\langle\Omega(X_\alpha,p):\alpha\in\delta\rangle$ is a decreasing
$\leRF$-chain.

To ensure that this chain does not have a lower bound one should make sure
that $p$~is not an accumulation point of a countable discrete subset
of the intersection.
Having a one-point intersection is certainly sufficient for this.
In~\cite{MR863908} Van Douwen showed that it is possible to have a chain
of length~$\cee$ with a one-point intersection.
In~\cite{MR633575} and~\cite{MR1007490} we find constructions of decreasing
$\leRF$-chains of type~$\omega$ and of type~$\mu$ for uncountable~$\mu<\cee$
respectively.
The latter two constructions provide a point in the intersection of
a suitable chain of copies of~$\betaN$ that is not an accumulation point of
a countable discrete subset of that intersection.

We want to know when in these cases the intersection can be made to be a
one-point set.
\end{question}

\begin{question}
If $\kappa\le\cee$ has uncountable cofinality and if
$\left\langle X_\alpha:\alpha<\kappa\right\rangle$ is a strictly
decreasing sequence of copies of~$\betaN$ with intersection~$K$,
is there a point~$p$ in~$K$ that is not an accumulation point of any
countable discrete subset of~$K$?
\comments
This is related to Question~\ref{q.one-point}: the chains of copies of~$\betaN$
in the positive results were chosen with care.
We want to know if that care is necessary.
\end{question}

\begin{question} What are the possible lengths of unbounded $\RF$-chains%
\pri{RF-chain@$\RF$-chain!unbounded}?
\comments
Since every point has at most~$\cee$ predecessors every chain has cardinality
at most~$\cee^+$.
In~\cite{MR633575} we find a point with exactly $\aleph_0$~many predecessors,
with the order type of the set of negative integers.

Every unbounded chain will have cardinality at least~$\cee$
(this follows from results in~\cite{MR277371}*{Theorem~2.9} or~\cite{MR698391}),
so the cardinality of an unbounded chain is equal to either~$\cee$ or~$\cee^+$.
In~\cite{MR730151} and~\cite{MR817833} Butkovi\v{c}ov\'a constructed unbounded
chains of order-type~$\cee^+$ and~$\omega_1$ respectively.

What other order-type are possible?
Can one prove that a chain or order-type~$\cee$ (or its cofinality)
exists, irrespective of~$\CH$?
\end{question}

\begin{question}
Is every finite partial order embeddable in the Rudin-Keisler order?
\comments
See MathOverFlow \texttt{https://mathoverflow.net/questions/375365}.
To get a positive answer it suffices to embed every finite power set
into this order.
It is relatively easy to adapt the construction of two incomparable
ultrafilters to yield an embedding of the power set of~$\{0,1\}$
(see~\cite{ruitje}),
but an embedding of the power set of~$\{0,1,2\}$ already poses unexpected
difficulties.\par
The analogous question for the Rudin-Frol\'ik order has an easier answer.
This order is tree-like in that the predecessors of a point are linearly
ordered, and every point has $2^\cee$ many successors.
This implies that every finite rooted tree, and only those among the finite
partial orders, can be embedded into this order.
\end{question}

\section{Uncountable Cardinals}

\begin{question}\label{q.small.chi}
Is there consistently an uncountable cardinal $\kappa$ with
$p\in U(\kappa)$ such that~$\chi(p)<2^\kappa$%
\pri{ultrafilter!of small character!on uncountable cardinals}?
\comments
It is well known that if $\kappa$~is an infinite cardinal then there
are $2^{2^\kappa}$~many uniform ultrafilters on~$\kappa$ with
character equal to~$2^\kappa$, see~\cites{MR1503375,MR0001454}.

It is also well known, and referred to in other questions, that it
is consistent that there are ultrafilters on~$\N$ of character less than~$\cee$.

Of course the Generalized Continuum Hypothesis implies that every uniform
ultrafilter on every~$\kappa$ has character~$2^\kappa$, but we are not aware
of any consistency result the other way for uncountable cardinals.

We formulate two special cases of our question:
\begin{itemize}
\item Is it consistent to have a uniform ultrafilter on~$\omega_1$ of
      character~$\aleph_2$ (with $\aleph_2<2^{\aleph_1}$ of course)?
\item Is it consistent to have a measurable cardinal $\kappa$ with
      a $p\in U(\kappa)$ such that $\chi(p)<2^\kappa$?
\end{itemize}
The first question simply looks at the smallest possible case and the second
question asks, implicitly, if having a uniform ultrafilter of small
character is actually a large-cardinal property of~$\aleph_0$.

There has been recent activity in this area;
the paper~\cite{MR3832086} deals with the character spectrum of uncountable
cardinals of countable cofinality, and in~\cite{MR4081063} one finds models
with $\you_\kappa<2^\kappa$ for $\kappa=\cee$ and for $\kappa=\aleph_{\omega+1}$.
These results use large cardinals in the ground model: the spectrum result
uses a supercompact and many measurables; the results for~$\cee$
and~$\aleph_{\omega+1}$ use a measurable and supercompact cardinal respectively.
\end{question}

\begin{question}
Is it consistent to have cardinals $\kappa<\lambda$ with points
$p\in U(\kappa)$ and $q\in U(\lambda)$ such that $\chi(p)>\chi(q)$%
\pri{ultrafilter!of small character!on uncountable cardinals}?
\comments
This is a follow-up question to Question~\ref{q.small.chi}:
if uniform ultrafilters of small character are at all possible, how
much variation can we achieve among various cardinals?
\end{question}

\begin{question}
If $\kappa$ is regular and uncountable, $\calF$ is a countably
complete uniform filter on $\kappa$ then what is the cardinality
of the closed set $U_\calF=\{u\in U(\kappa):\calF\subseteq u\}$?
\comments
In case $\kappa$ is measurable one can use a measure ultrafilter to create
filters~$\calF$ such that $U_\calF$ is finite or a copy of~$\beta\lambda$,
for any $\lambda<\kappa$.

For other cardinals the set $U_\calF$ will always be at least infinite and
given the nature of~$\beta\kappa$ the cardinality will be closely related
to numbers of the form~$2^{2^\lambda}$ for $\lambda\le\kappa$.

For the closed unbounded filter the answer is~$2^{2^\kappa}$: using
a family of $\kappa$~many pairwise disjoint stationary sets and an independent
family on~$\kappa$ of cardinality~$2^\kappa$ one can produce a map from~$U_\calF$
onto the Cantor cube of weight~$2^\kappa$.
\end{question}

\begin{question}
Assume that $\kappa$ is regular, that $\kappa\subseteq X\subseteq\beta\kappa$
is such that $[X]_{<\kappa}=X$ and $\beta_X\kappa=X$.
Now if $Y$ is a closed subspace of a power of $X$, is then also
$X$ a closed subspace of a power of $Y$?
\comments
Some notation: $[X]_{<\kappa}$
denotes $\bigcup\{\cl B:B\in[X]^{<\kappa}\}$, and
if $\kappa\subseteq X\subseteq\beta\kappa$ then $\beta_X\kappa$~is
the maximal subset of~$\beta\kappa$ such that every function from~$\kappa$
to~$X$ has a continuous extension from~$\beta_X\kappa$ to~$X$.
\end{question}

\begin{question}
Are there $\kappa$ and $p\in U(\kappa)$ such that $|\R_p|>|\R_p/\equiv|=\cee$%
\pri{cardinality!of ultrapower}\pri{ultrapower!cardinality of}?
\comments
Here $\R_p$ denotes the ultrapower of~$\R$ by the ultrafilter~$p$.
The relation~$\equiv$ is that of Archimedean equivalence:
$a\equiv b$ means that there is an~$n\in\N$ such that both $\abs{a}<\abs{nb}$
and $\abs{b}<\abs{na}$.
\end{question}

\begin{question}
Is there a $C^*$-embedded bi-Bernstein set\pri{bi-Bernstein set}
in $U(\omega_1)$\pri{Uomega1@$U(\omega_1)$}?
\comments
\end{question}

\begin{question}
Are there open sets $G_1$ and $G_2$ in $U(\omega_1)$%
\pri{open set!in $U(\omega_1)$} such that
$\cl{G_1}\cap\cl{G_2}$ consists of exactly one point?
\comments
\end{question}

\end{document}